# A Generalization of Obreshkoff-Ehrlich Method for Multiple Roots of Algebraic, Trigonometric and Exponential Equations

*A. I. Iliev*

In this paper methods for simultaneous finding all roots of generalized polynomials are developed. These methods are related to the case when the roots are multiple. They possess cubic rate of convergence and they are as labour-consuming as the known methods related to the case of polynomials with simple roots only.

**Introduction.** After 1960 the question of simultaneous finding all roots (SFAR) of polynomials became very actual and it is considered by many authors. The reason of this interest is the better behaviour of the methods for SFAR with respect to the methods for individual search of the roots. Also these methods are very convenient for application on the computers with parallel processors. Methods for SFAR have a wider region of convergence and they are more stable. In several survey publications [1,2,3] this question is considered in details. The first methods for SFAR are related to the case when the roots are simple. The well-known method of Dochev [4] is for SFAR of algebraic polynomial with real and simple roots. The developments of this same method for the case of nonalgebraic polynomials ( trigonometric, exponential and generalized ) are performed in [5,6,7]. The classical method of Obreshkoff-Ehrlich [8] possessing cubic rate of convergence is also generalized [9]. Using the approach basing on the divided differences with multiple knots Semerdzhiev [10] generalized the method of Dochev to the case when the roots have arbitrary, but given multiplicities. The same question for the case of trigonometric and exponential polynomials is solved in [11,3]. The new methods preserve their quadratic rate of convergence. The method of Obreshkoff-Ehrlich is also generalized to the most general case [12,3] of polynomials upon some Chebyshev system, having multiple roots with given multiplicities. The rate of convergence is cubic but, unfortunately, this generalization requires at each iteration to calculate determinants which is a labour-consuming operation.

In this paper we develop a new method which is a generalization of the Obreshkoff-Ehrlich method for the cases of algebraic, trigonometric and exponential polynomials. This method has a cubic rate of convergence. It is efficient from the computational point of view and can be used for SFAR if the roots have known multiplicities. This new method in spite of the arbitrariness of multiplicities is of the same complexity as the methods for SFAR of simple roots. We do not use divided differences with multiple knots and this fact does not lead to calculation of derivatives of the given polynomial of higher order, but only of first ones. The results of this paper are published in shortened form as a preliminary communication in [13].

**Algebraic polynomials.** Let the algebraic polynomial

(1) $\quad A_n(x) = x^n + a_1 x^{n-1} + \ldots + a_n$

be given and let $x_1, x_2, \ldots, x_m$ be its roots with given multiplicities $\alpha_1, \alpha_2, \ldots, \alpha_m$ respectively $(\alpha_1 + \alpha_2 + \ldots + \alpha_m = n)$. For SFAR of (1) the Ehrlich formula

(2) $\quad x_i^{[k+1]} = x_i^{[k]} - A_n\left(x_i^{[k]}\right)\left[ A_n'\left(x_i^{[k]}\right) - A_n\left(x_i^{[k]}\right) \sum_{j=1, j\neq i}^{n} \left(x_i^{[k]} - x_j^{[k]}\right)^{-1} \right]^{-1}$

$i = \overline{1,n}, k = 0,1,2,\ldots$

is well known. Formula (2) can be written in the form

(3) $\quad x_i^{[k+1]} = x_i^{[k]} - A_n\left(x_i^{[k]}\right)\left[ A_n'\left(x_i^{[k]}\right) - A_n\left(x_i^{[k]}\right) Q''^{[k]}\left(x_i^{[k]}\right)\left[ 2 Q'^{[k]}\left(x_i^{[k]}\right) \right]^{-1} \right]^{-1}$

$i = \overline{1,n}, k = 0,1,2,\ldots$

where

(4) $\quad Q^{[k]}(x) = \prod_{j=1}^{n} \left(x - x_j^{[k]}\right).$

We define

(5) $\quad x_i^{[k+1]} = x_i^{[k]} - \alpha_i A_n\left(x_i^{[k]}\right)\left[ A_n'\left(x_i^{[k]}\right) - A_n\left(x_i^{[k]}\right) Q_i'^{[k]}\left(x_i^{[k]}\right) / Q_i^{[k]}\left(x_i^{[k]}\right) \right]^{-1}$

$i = \overline{1,m}, k = 0,1,2\ldots$

where

(6) $\quad Q_i^{[k]}(x) = \prod_{j=1, j\neq i}^{m} \left(x - x_j^{[k]}\right)^{\alpha_j}.$

**Theorem 1.** *Let $q$, $c$ and $d \stackrel{\text{def}}{=} \min_{i \neq j} |x_i - x_j|$ be real constants such that the following inequalities*

(7) $\quad \begin{aligned} & 1 > q > 0, \quad c > 0, \quad d - 2c > 0, \\ & 0 < c^2(n - 3\alpha_i) + (n + (3d-1)\alpha_i)c < d^2 \alpha_i, \quad i = \overline{1,m} \end{aligned}$

*be satisfied. If the initial approximations $x_1^{[0]}, \ldots, x_m^{[0]}$ to the exact roots $x_1, \ldots, x_m$ of (1) are chosen so that the inequalities $\left|x_i^{[0]} - x_i\right| \leq cq$, $i = \overline{1,m}$ hold true then for every natural $k$ the inequalities*

(8) $\quad \left|x_i^{[k]} - x_i\right| \leq c q^{3^k}$, $\quad i = \overline{1,m}$

*also hold true.*

Proof.

We prove the theorem 1 by means of induction with respect to the number of the iteration. From the assumptions of the theorem we have that (8) are fulfilled for $k=0$. Suppose that (8) hold true for some $k>0$.

From (5) we obtain

(9) $\quad x_i^{[k+1]} - x_i = x_i^{[k]} - x_i - \alpha_i \left[ A_n'(x_i^{[k]}) / A_n(x_i^{[k]}) - Q_i'^{[k]}(x_i^{[k]}) / Q_i^{[k]}(x_i^{[k]}) \right]^{-1} \quad , i = \overline{1,m}$.

Using the representation $Q_i'^{[k]}(x_i^{[k]}) / Q_i^{[k]}(x_i^{[k]}) = \sum_{j=1, j \neq i}^{m} \alpha_j / (x - x_j^{[k]}) \quad , i = \overline{1,m}$

from (9) we receive

$$x_i^{[k+1]} - x_i = x_i^{[k]} - x_i - \alpha_i \left[ \sum_{j=1}^{m} \alpha_j / (x_i^{[k]} - x_j) - \sum_{j=1, j \neq i}^{m} \alpha_j / (x_i^{[k]} - x_j^{[k]}) \right]^{-1} \quad , i = \overline{1,m}.$$

Further transformations lead to

$$x_i^{[k+1]} - x_i = (x_i^{[k]} - x_i) \left[ 1 - \alpha_i \left[ \alpha_i + (x_i^{[k]} - x_i) \sum_{j=1, j \neq i}^{m} \alpha_j (x_j - x_j^{[k]}) \left[ (x_i^{[k]} - x_j)(x_i^{[k]} - x_j^{[k]}) \right]^{-1} \right]^{-1} \right] =$$

(10) $\quad (x_i^{[k]} - x_i)^2 \left[ \alpha_i + (x_i^{[k]} - x_i) \sum_{j=1, j \neq i}^{m} \alpha_j (x_j - x_j^{[k]}) \left[ (x_i^{[k]} - x_j)(x_i^{[k]} - x_j^{[k]}) \right]^{-1} \right]^{-1} \times$

$$\sum_{j=1, j \neq i}^{m} \alpha_j (x_j - x_j^{[k]}) \left[ (x_i^{[k]} - x_j)(x_i^{[k]} - x_j^{[k]}) \right]^{-1} \quad , i = \overline{1,m}.$$

Obviously we have

(11) $\quad \left| x_i^{[k]} - x_j \right| \geq \left| x_i - x_j \right| - \left| x_i - x_i^{[k]} \right| \geq d - cq^{3^k} > d - c \quad , i, j = \overline{1,m} \quad , i \neq j$.

On the other hand

(12) $\quad \left| x_i^{[k]} - x_j^{[k]} \right| \geq \left| x_i^{[k]} - x_j \right| - \left| x_j - x_j^{[k]} \right| \geq d - 2cq^{3^k} > d - 2c \quad , i, j = \overline{1,m} \quad , i \neq j$.

Using (10)-(12) we find as a final result

$$\left| x_i^{[k+1]} - x_i \right| \leq \left| x_i^{[k]} - x_i \right|^2 \left[ \alpha_i - \left| x_i^{[k]} - x_i \right| \sum_{j=1, j \neq i}^{m} \alpha_j \left| x_j - x_j^{[k]} \right| \left[ \left| x_i^{[k]} - x_j \right| \left| x_i^{[k]} - x_j^{[k]} \right| \right]^{-1} \right]^{-1} \times$$

$$\sum_{j=1, j \neq i}^{m} \alpha_j \left| x_j - x_j^{[k]} \right| \left[ \left| x_i^{[k]} - x_j \right| \left| x_i^{[k]} - x_j^{[k]} \right| \right]^{-1} \leq$$

$$c^3 (q^{3^k})^3 \left[ \alpha_i - c \sum_{j=1, j \neq i}^{m} \alpha_j [(d-c)(d-2c)]^{-1} \right]^{-1} \sum_{j=1, j \neq i}^{m} \alpha_j [(d-c)(d-2c)]^{-1} =$$

$$cq^{3^{k+1}} c^2 \left[ \alpha_i - c(n - \alpha_i)[(d-c)(d-2c)]^{-1} \right]^{-1} (n - \alpha_i)[(d-c)(d-2c)]^{-1} < cq^{3^{k+1}}$$

$i = \overline{1,m}$

which proves the theorem completely.

**Proposition 1.** *In the case when* $\alpha_1 = \alpha_2 = ... = \alpha_m = 1$ *then the correlations*

(13) $\quad Q''^{[k]}(x_i^{[k]}) / Q'^{[k]}(x_i^{[k]}) = 2 Q_i'^{[k]}(x_i^{[k]}) / Q_i^{[k]}(x_i^{[k]}) \quad , i = \overline{1,m}$

*hold true.*

Proof. In this case (6) reduces itself to (4) and we have
(14) $\quad Q^{[k]}(x) = Q_i^{[k]}(x)(x - x_i^{[k]}) \quad , i = \overline{1,m}.$

Differentiating (14) we obtain
(15)
$$Q'^{[k]}(x) = Q_i'^{[k]}(x)(x - x_i^{[k]}) + Q_i^{[k]}(x)$$
$$Q''^{[k]}(x) = Q_i''^{[k]}(x)(x - x_i^{[k]}) + 2Q_i'^{[k]}(x) \quad , i = \overline{1,m}.$$

From (15) we receive (13).

Proposition 1 shows that method (5) coincides with the method (3) and consequently with the method of Ehrlich in the case when $\alpha_i = 1$ , $i = \overline{1,m}$.

Example 1. For the equation $A_6(x) = (x+2)^2(x-1)(x-3)^3 = 0$ at the initial approximations $x_1^{[0]} = -3$ , $x_2^{[0]} = 0.1$ and $x_3^{[0]} = 4$ using the formula (5) we receive the roots with 18 decimal digits after only 4 iterations.

| k | $x_1^{[k]}$ | $x_2^{[k]}$ | $x_3^{[k]}$ |
|---|---|---|---|
| 0 | -3.00000000000000000 | 0.10000000000000000 | 4.00000000000000000 |
| 1 | -1.99942363112391931 | 1.03532819268537456 | 3.03985932004689332 |
| 2 | -2.00000000143304088 | 0.999961906975802837 | 2.99999539984403290 |
| 3 | -2.00000000000000000 | 1.00000000000000501 | 3.00000000000000007 |
| 4 | -2.00000000000000000 | 1.00000000000000000 | 3.00000000000000000 |

**Trigonometric polynomials.** For the trigonometric polynomial $T_n(x) = a_0/2 + \sum_{l=1}^{n}(a_l \cos(lx) + b_l \sin(lx))$ we suppose that at least one of the leading coefficients $a_n$ and $b_n$ is not zero and that it has real roots $x_1,...,x_m$ with given multiplicities $\alpha_1, \alpha_2,...,\alpha_m$ $(\alpha_1 + \alpha_2 + ... + \alpha_m = 2n)$. Analogously to (5) we can use the iteration method

(16)
$$x_i^{[k+1]} = x_i^{[k]} - \alpha_i T_n(x_i^{[k]})\left[T_n'(x_i^{[k]}) - T_n(x_i^{[k]})Q_i'^{[k]}(x_i^{[k]})/Q_i^{[k]}(x_i^{[k]})\right]^{-1}$$
$$i = \overline{1,m} \quad , \quad k = 0,1,2,...$$

where
$$Q_i^{[k]}(x) = \prod_{j \neq i, j=1}^{m} \sin^{\alpha_j}\left((x - x_j^{[k]})/2\right).$$

The formula (16) at $\alpha_1 = \alpha_2 = ... = \alpha_m = 1$ coincides with the analogue of Obreshkoff-Ehrlich formula [9] for trigonometric polynomials.

**Theorem 2.** *Let us denote* $d \stackrel{\text{def}}{=} \min_{i \neq j}|x_i - x_j|$. *Let* $c$, $q$ *and* $\xi$ *be positive real numbers so that* $q < 1$, $2c < \xi$, $d - 2c > 0$ *and* $\max_{i \neq j}|x_i - x_j| < 2\pi - 2\xi$. *Denote the expression* $\min\{|\sin \xi/2|, |\sin(d/2 - c)|\}$ *by* A. *If*

$c^2(4n + \alpha_i(9A^2/8 - 2)) < A^2\alpha_i$, $i = \overline{1,m}$ and initial approximations $x_i^{[0]}$, $i = \overline{1,m}$ are chosen so that $|x_i^{[0]} - x_i| \leq cq$, $i = \overline{1,m}$ then for every natural $k$ the inequalities $|x_i^{[k]} - x_i| \leq cq^{3^k}$, $i = \overline{1,m}$ also hold true.

Proof. We divide the numerator and denominator of second summand in the right side of (16) by $T_n(x_i^{[k]})$ and obtain

(17) $\quad x_i^{[k+1]} - x_i = x_i^{[k]} - x_i - \alpha_i \left[ T_n'(x_i^{[k]}) / T_n(x_i^{[k]}) - Q_i'^{[k]}(x_i^{[k]}) / Q_i^{[k]}(x_i^{[k]}) \right]^{-1}$, $i = \overline{1,m}$.

On the other hand we have

(18) $\quad T_n'(x_i^{[k]}) / T_n(x_i^{[k]}) = 2^{-1} \sum_{j=1}^{m} \alpha_j \, \text{cotg}\left((x_i^{[k]} - x_j)/2\right)$, $i = \overline{1,m}$

and

(19) $\quad Q_i'^{[k]}(x_i^{[k]}) / Q_i^{[k]}(x_i^{[k]}) = 2^{-1} \sum_{j=1, j \neq i}^{m} \alpha_j \, \text{cotg}\left((x_i^{[k]} - x_j^{[k]})/2\right)$, $i = \overline{1,m}$.

Using (18) and (19) we transform (17) into the form

(20) $\quad x_i^{[k+1]} - x_i = x_i^{[k]} - x_i - 2\alpha_i \left[ \sum_{j=1}^{m} \alpha_j \, \text{cotg}\left((x_i^{[k]} - x_j)/2\right) - \sum_{j=1, j \neq i}^{m} \alpha_j \, \text{cotg}\left((x_i^{[k]} - x_j^{[k]})/2\right) \right]^{-1}$

$i = \overline{1,m}$.

If we multiply the numerator and denominator of the second summand in the right side of (20) with $\sin\left((x_i^{[k]} - x_i)/2\right)$ then we obtain

$x_i^{[k+1]} - x_i = x_i^{[k]} - x_i - 2\alpha_i \times$

$\left[ \alpha_i \cos\left((x_i^{[k]} - x_i)/2\right) + \sin\left((x_i^{[k]} - x_i)/2\right) \sum_{j=1, j \neq i}^{m} \alpha_j \left[ \text{cotg}\left((x_i^{[k]} - x_j)/2\right) - \text{cotg}\left((x_i^{[k]} - x_j^{[k]})/2\right) \right] \right]^{-1} \times$

$\sin\left((x_i^{[k]} - x_i)/2\right) = (x_i^{[k]} - x_i) \times$

$\left[ \alpha_i \left( \cos\left((x_i^{[k]} - x_i)/2\right) - 2(x_i^{[k]} - x_i)^{-1} \sin\left((x_i^{[k]} - x_i)/2\right) \right) + \sin\left((x_i^{[k]} - x_i)/2\right) \times \right.$

$\sum_{j=1, j \neq i}^{m} \alpha_j \left[ \text{cotg}\left((x_i^{[k]} - x_j)/2\right) - \text{cotg}\left((x_i^{[k]} - x_j^{[k]})/2\right) \right] \right] \times$

$\left[ \alpha_i \cos\left((x_i^{[k]} - x_i)/2\right) + \sin\left((x_i^{[k]} - x_i)/2\right) \sum_{j=1, j \neq i}^{m} \alpha_j \left[ \text{cotg}\left((x_i^{[k]} - x_j)/2\right) - \text{cotg}\left((x_i^{[k]} - x_j^{[k]})/2\right) \right] \right]^{-1}$, $i = \overline{1,m}$.

Further the difference $\text{cotg}\left((x_i^{[k]} - x_j)/2\right) - \text{cotg}\left((x_i^{[k]} - x_j^{[k]})/2\right)$ can be transformed as follows

$$\cotg\left(\left(x_i^{[k]} - x_j\right)/2\right) - \cotg\left(\left(x_i^{[k]} - x_j^{[k]}\right)/2\right) = \left[\sin\left(\left(x_i^{[k]} - x_j\right)/2\right)\sin\left(\left(x_i^{[k]} - x_j^{[k]}\right)/2\right)\right]^{-1} \times$$

$$\left[\cos\left(\left(x_i^{[k]} - x_j\right)/2\right)\sin\left(\left(x_i^{[k]} - x_j^{[k]}\right)/2\right) - \cos\left(\left(x_i^{[k]} - x_j^{[k]}\right)/2\right)\sin\left(\left(x_i^{[k]} - x_j\right)/2\right)\right] =$$

$$\sin\left(\left(x_j - x_j^{[k]}\right)/2\right)\left[\sin\left(\left(x_i^{[k]} - x_j\right)/2\right)\sin\left(\left(x_i^{[k]} - x_j^{[k]}\right)/2\right)\right]^{-1}, i,j = \overline{1,m}, i \neq j.$$

Consequently, for the deviation of $x_i^{[k+1]}$ from $x_i$ we receive the expression

(21) $x_i^{[k+1]} - x_i = \left[\alpha_i\left[\left(x_i^{[k]} - x_i\right)\cos\left(\left(x_i^{[k]} - x_i\right)/2\right) - 2\sin\left(\left(x_i^{[k]} - x_i\right)/2\right)\right] + \left(x_i^{[k]} - x_i\right)Y_i^{[k]}\left(x_i^{[k]}\right)\right] \times$

$$\left[\alpha_i \cos\left(\left(x_i^{[k]} - x_i\right)/2\right) + Y_i^{[k]}\left(x_i^{[k]}\right)\right]^{-1}$$

$$Y_i^{[k]}\left(x_i^{[k]}\right) = \sin\left(\left(x_i^{[k]} - x_i\right)/2\right) \sum_{j=1, j\neq i}^{m} \alpha_j \sin\left(\left(x_j - x_j^{[k]}\right)/2\right)\left[\sin\left(\left(x_i^{[k]} - x_j\right)/2\right)\sin\left(\left(x_i^{[k]} - x_j^{[k]}\right)/2\right)\right]^{-1}$$

$$i = \overline{1,m}.$$

In order to find an estimate for the expressions
$\left(x_i^{[k]} - x_i\right)\cos\left(\left(x_i^{[k]} - x_i\right)/2\right) - 2\sin\left(\left(x_i^{[k]} - x_i\right)/2\right), i = \overline{1,m}$, we consider the auxiliary function $F(x) = x\cos(x/2) - 2\sin(x/2)$ and its Mak-Laurent expansion until the remainder term with third derivative of $F(x)$. In this way we obtain

$$F\left(x_i^{[k]} - x_i\right) = \left(-2^{-1}\cos\left(\zeta_i^{[k]}/2\right) + \left(\zeta_i^{[k]}/8\right)\sin\left(\zeta_i^{[k]}/2\right)\right)\left(x_i^{[k]} - x_i\right)^3/6$$

$$\left(\zeta_i^{[k]} = \theta_i^{[k]}\left(x_i^{[k]} - x_i\right), 0 < \theta_i^{[k]} < 1, i = \overline{1,m}\right)$$

and, therefore, the following estimate

$$\left|F\left(x_i^{[k]} - x_i\right)\right| \leq (1/12 + 2\pi/48)\left|x_i^{[k]} - x_i\right|^3 \leq \left|x_i^{[k]} - x_i\right|^3/4 \leq \left|x_i^{[k]} - x_i\right|^3, i = \overline{1,m}.$$

On the other hand the inequalities (11) and (12) hold true. Because of the fact that all the roots are in an interval with a length $2\pi$ i.e. $\left|x_i - x_j\right| < 2\pi, i,j = \overline{1,m}, i \neq j$ it exists a positive number $\xi$ such that $\left|x_i - x_j\right| < 2\pi - 2\xi, i,j = \overline{1,m}, i \neq j$. We now receive the inequalities

$$\left|x_i^{[k]} - x_j^{[k]}\right| \leq \left|x_i^{[k]} - x_i\right| + \left|x_j^{[k]} - x_j\right| + \left|x_i - x_j\right| < 2\pi - 2\xi + 2cq^{3^k}, i,j = \overline{1,m}, i \neq j$$

$$\left|x_i^{[k]} - x_j\right| \leq \left|x_i^{[k]} - x_i\right| + \left|x_i - x_j\right|, i,j = \overline{1,m}, i \neq j.$$

From the suppositions of the theorem it follows that $2c < \xi$. Then

$$\left|x_i^{[k]} - x_j^{[k]}\right| < 2\pi - \xi, i,j = \overline{1,m}, i \neq j$$

$$\left|x_i^{[k]} - x_j\right| < 2\pi - 2\xi + \xi/2 < 2\pi - \xi, i,j = \overline{1,m}, i \neq j.$$

Consequently $d/2 - c < \left|x_i^{[k]} - x_j\right|/2 < \pi - \xi/2, i,j = \overline{1,m}, i \neq j$ and
$d/2 - c < \left|x_i^{[k]} - x_j^{[k]}\right|/2 < \pi - \xi/2, i,j = \overline{1,m}, i \neq j$. It is easy to find that both expressions $\left|\sin\left(\left(x_i^{[k]} - x_j^{[k]}\right)/2\right)\right|$ and $\left|\sin\left(\left(x_i^{[k]} - x_j\right)/2\right)\right|$ are greater than $A$,

$A \overset{\text{def}}{=} \min\{|\sin(\xi/2)|, |\sin(d/2-c)|\}$. From (21) we estimate the absolute value of $x_i^{[k+1]} - x_i, i = \overline{1,m}$ i.e.

$$|x_i^{[k+1]} - x_i| \leq \left[\alpha_i |x_i^{[k]} - x_i|^3 + |x_i^{[k]} - x_i| Z_i^{[k]}(x_i^{[k]})\right]\left[\alpha_i |\cos((x_i^{[k]} - x_i)/2)| - Z_i^{[k]}(x_i^{[k]})\right]^{-1}$$

(22) $Z_i^{[k]}(x_i^{[k]}) = |\sin((x_i^{[k]} - x_i)/2)| \sum_{j=1, j \neq i}^{m} \alpha_j |\sin((x_j - x_j^{[k]})/2)| \left[|\sin((x_i^{[k]} - x_j)/2)| |\sin((x_i^{[k]} - x_j^{[k]})/2)|\right]^{-1}$

$i = \overline{1,m}$.

Because of the presentation $\sin((x_i^{[k]} - x_i)/2) = ((x_i^{[k]} - x_i)/2)\cos\varsigma_i^{[k]}, i = \overline{1,m}$, where $\varsigma_i^{[k]} = \theta_i^{[k]}((x_i^{[k]} - x_i)/2), 0 < \theta_i^{[k]} < 1, i = \overline{1,m}$ the estimates $|\sin((x_i^{[k]} - x_i)/2)| \leq |x_i^{[k]} - x_i|/2, i = \overline{1,m}$ hold true. Then from (22) we obtain

$$|x_i^{[k+1]} - x_i| \leq c^3(q^{3^k})^3 \left[\alpha_i + A^{-2} \sum_{j=1, j\neq i}^{m} \alpha_j\right]\left[\alpha_i |\cos((x_i^{[k]} - x_i)/2)| - (c/A)^2 \sum_{j=1, j\neq i}^{m} \alpha_j\right]^{-1}, i = \overline{1,m}.$$

From the inequality $|(x_i^{[k]} - x_i)/2| < c/2, i = \overline{1,m}$ and the presentation $\cos((x_i^{[k]} - x_i)/2) = 1 - (1/8)(x_i^{[k]} - x_i)^2 \cos\zeta_i^{[k]}, i = \overline{1,m}$ it follows that $|\cos((x_i^{[k]} - x_i)/2)| > 1 - c^2/8, i = \overline{1,m}$ for sufficiently small $c$. Finally we receive

$$|x_i^{[k+1]} - x_i| \leq cq^{3^{k+1}} c^2 \left[\alpha_i + (2n - \alpha_i)/A^2\right]\left[\alpha_i(1 - c^2/8) - (2n - \alpha_i)(c/A)^2\right]^{-1} < cq^{3^{k+1}}, i = \overline{1,m}$$

for a small enough $c$. Thus the theorem is proved completely.

Example 2. For the trigonometric polynomial
$$T_3(x) = \sin^3((x-1)/2)\sin^2((x-2)/2)\sin((x-2.5)/2)$$
at initial approximations $x_1^{[0]} = 0.2$, $x_2^{[0]} = 1.7$ and $x_3^{[0]} = 3$ we reach the roots of $T_3(x)$ with an accuracy of 18 digits at the $5^{th}$ iteration.

| k | $x_1^{[k]}$ | $x_2^{[k]}$ | $x_3^{[k]}$ |
|---|---|---|---|
| 0 | 0.200000000000000000 | 1.700000000000000000 | 3.000000000000000000 |
| 1 | 1.08093197781206681 | 2.13081574593339511 | 2.68530050098035859 |
| 2 | 0.999087999636487434 | 1.98917328088624173 | 2.46587439388854078 |
| 3 | 1.00000001182848523 | 2.00000867262537340 | 2.50012119040535689 |
| 4 | 1.00000000000000000 | 1.99999999999998133 | 2.49999999999881136 |
| 5 | 1.00000000000000000 | 2.00000000000000000 | 2.50000000000000000 |

**Exponential polynomials.** Let us now consider the polynomial

(23) $E_n(x) = a_0/2 + \sum_{l=1}^{n}(a_l \operatorname{ch}(lx) + b_l \operatorname{sh}(lx)) = a_0/2 + \sum_{l=1}^{n}(a'_l \ell^{lx} + b'_l \ell^{-lx})$.

We suppose that at least one of the leading coefficients $a_n$ or $b_n$ is not zero and that $E_n(x)$ has real roots $x_1, x_2, \ldots, x_m$ with known multiplicities

$\alpha_1, \alpha_2, \ldots, \alpha_m$ $(\alpha_1 + \alpha_2 + \ldots + \alpha_m = 2n)$ correspondingly. The roots of (23) can be refined simultaneously with the help of the computational scheme

$$(24) \quad x_i^{[k+1]} = x_i^{[k]} - \alpha_i E_n(x_i^{[k]}) \left[ E_n'(x_i^{[k]}) - E_n(x_i^{[k]}) Q_i'^{[k]}(x_i^{[k]}) / Q_i^{[k]}(x_i^{[k]}) \right]^{-1}$$

$$i = \overline{1, m}, \quad k = 0, 1, 2, \ldots$$

where

$$Q_i^{[k]}(x) = \prod_{j \neq i, j=1}^{m} \operatorname{sh}^{\alpha_j}\left( (x - x_j^{[k]})/2 \right).$$

**Theorem 3.** *Denote $\min_{i \neq j} |x_i - x_j|$ by $d$. Let $q$ and $c$ be real numbers such that $1 > q > 0$, $c > 0$, $d - 2c > 0$, $\left(2^{-1} \operatorname{ch}(c/2) + 8^{-1} c |\operatorname{sh}(c/2)|\right) < 6$, $\operatorname{ch}(c/2) < 2$ and $c^2(4n + (S^2 - 2)\alpha_i) < S^2 \alpha_i$, $i = \overline{1, m}$, where the expression $\operatorname{sh}((d-2c)/2)$ is indicated by $S$. If the initial approximations $x_i^{[0]}$, $i = \overline{1, m}$ are taken such that $|x_i^{[0]} - x_i| \leq cq$, $i = \overline{1, m}$ then for every $k \in N$ the inequalities $|x_i^{[k]} - x_i| \leq cq^{3^k}$, $i = \overline{1, m}$ hold true.*

The proof of theorem 3 can be carried out by the similar manner as the theorem 2 with corresponding changes which are connected with the properties of the hyperbolic functions.

Example 3. The iteration method (24) was applied for SFAR of the exponential polynomial $E_2(x) = \operatorname{sh}^2((x+2)/2) \operatorname{sh}^2((x-3)/2)$. Using initial approximations $x_1^{[0]} = -1$ and $x_2^{[0]} = 4$ by the formula (24) we receive the roots with 18 decimal digits after only 4 iterations.

| $k$ | $x_1^{[k]}$ | $x_2^{[k]}$ |
| --- | --- | --- |
| 0 | -1.00000000000000000 | 4.00000000000000000 |
| 1 | -1.93448948248966207 | 3.07207901269406155 |
| 2 | -1.99997875689833755 | 3.00002895806496640 |
| 3 | -1.99999999999999929 | 3.00000000000000190 |
| 4 | -2.00000000000000000 | 3.00000000000000000 |

*Acknowledgements.* The author is deeply grateful to Khristo Semerdzhiev for formulating the problem and for many helpful discussions and advice during this investigation.

*University of Plovdiv* http://www.pu.acad.bg
*Faculty of Mathematics and Informatics* http://www.fmi.pu.acad.bg
*Tzar Assen Str.,* 24
4000 *Plovdiv*
*BULGARIA*
*e-mail:* aii@pu.acad.bg
*URL:* http://anton.iliev.tripod.com